&amslplain
\documentstyle[amssymb,11pt]{amsart}
\newcommand{\half}{\frac12}
\newcommand{\Z}{{\Bbb Z}}

\newcommand{\C}{{\Bbb C}}

\newcommand{\Ref}[1]{{$($\ref{#1}$)$}}
\newcommand{\bean}{\begin{eqnarray}}
\newcommand{\eean}{\end{eqnarray}}
\newcommand{\be}{\begin{displaymath}}
\newcommand{\ee}{\end{displaymath}}
\newcommand{\bea}{\begin{eqnarray*}}
\newcommand{\eea}{\end{eqnarray*}}
\newcommand{\g}{{{\frak g}\,}}

\newcommand{\h}{{{\frak h\,}}}
\newcommand{\Id}{{\operatorname{Id}}}

\newcommand{\vs}{\vspace{1.5\baselineskip}}
\newenvironment{prf}{\noindent{\it Proof\/}:}{$\;\Box$
\par\vs}

\newtheorem%
{thm}{Theorem}
\newtheorem%
{proposition}[thm]{Proposition}
\newtheorem%
{lemma}[thm]{Lemma}
\newtheorem%
{lemmadef}[thm]{Lemma-Definition}
\newtheorem%
{corollary}[thm]{Corollary}
\newtheorem%
{conjecture}[thm]{Conjecture}
\newcommand{\End}{{\operatorname{End}}}

\newcommand{\tr}{{\operatorname{tr}}}
\newcommand{\Fun}{{\operatorname{Fun}}}
\title[Algebraic Bethe ansatz]
{Algebraic Bethe ansatz for the elliptic quantum group 
$E_{\tau,\eta}(sl_2)$}
\thanks{The authors were supported in part by NSF grants
DMS-9400841 and DMS-9501290}
\author{G. Felder and A. Varchenko}
\date{May 1996}
\address{Department of Mathematics,
University of North Carolina at Chapel Hill,
Chapel Hill, NC 27599-3250}
\begin{document}
\maketitle
\begin{abstract}
To each representation of the elliptic quantum group
$E_{\tau,\eta}(sl_2)$ is associated a family of commuting transfer
matrices. We give common eigenvectors by a version of the 
algebraic Bethe
ansatz method. Special cases of this construction give
eigenvectors for IRF models, for the eight-vertex model and
for the two-body Ruijsenaars operator. The latter is a
$q$-deformation of Hermite's solution of the Lam\'e equation.
\end{abstract}
\section{Introduction}
The Bethe ansatz is a method to construct common eigenvectors
of commuting families of operators (transfer matrices) occurring
in two-dimensional models of statistical mechanics. 
Faddeev and Takhtadzhan \cite{TF} reformulated the problem into
a question of representation theory: commuting families
of transfer matrices are associated to representations of
 certain algebras with
quadratic relations (now called quantum groups). Eigenvectors
are constructed by properly acting with algebra elements on ``highest
weight vectors''. In this form, the Bethe ansatz is
called algebraic Bethe ansatz.

Whereas this construction has been very successful in 
rational and trigonometric integrable models, its extension
to elliptic models has been problematic, although the
Bethe ansatz, in the case of the eight-vertex model,
is known since Baxter's work \cite{Ba}. For elliptic models 
associated to $sl_N$, an algebra with quadratic relations
has been introduced by Sklyanin \cite{S}, but the notion
of highest weight vector is not defined for its representations,
which makes a direct application of the algebraic Bethe ansatz
impossible. 

Recently, a definition of elliptic quantum groups $E_{\tau,\eta}(\g)$
associated
to any simple classical Lie algebra $\g$ was given \cite{Fe}. 
It is related to
a $q$-deformation of the Knizhnik--Zamolodchikov--Bernard
equation on tori. The representation theory of $E_{\tau,\eta}(sl_2)$
was described in \cite{FV}. 

In this paper we describe the algebraic Bethe ansatz for
$E_{\tau,\eta}(sl_2)$. The construction is very close to 
the construction done in the trigonometric case in \cite{TF}. The
main difference is that transfer matrices act on spaces of vector-valued
functions
rather than on finite dimensional vector spaces.

The basic results in this paper are Theorems \ref{three} and
\ref{four}. We present them in Sect.\ \ref{seconds} after 
a summary of the representation theory of $E_{\tau,\eta}(sl_2)$
in Sect.\ \ref{firsts}. In Sect.\ \ref{thirds} we formulate our
result in the discrete case, which gives a construction
of eigenfunctions for interaction-round-a-face models (Sect.\ 
\ref{fourths}).

In Sect.\ \ref{fifths} we introduce a version of Baxter's
vertex-IRF transformation. With the help of this, we show
how to obtain eigenvectors of the eight-vertex model transfer
matrix
from our eigenvectors. From this point of view, our construction
is very reminiscent to Baxter's transformation of the eight-vertex
model into an inhomogeneous six-vertex model \cite{Ba}.
Also, this vertex-IRF transformation shows a very close
relation of $E_{\tau,\eta}$ with 
Sklyanin's algebra in the $sl_2$ case. Note however that 
Sklyanin's algebra can only be defined for $sl_N$, whereas
the elliptic quantum group exists (at least) for all classical
simple Lie algebras.

Another class of problems in which the Bethe ansatz has been
applied successfully is the class of Calogero--Moser--Sutherland
quantum many body problems on the line with elliptic potentials,
see \cite{FVIMRN}, \cite{FVComp}. In the case of two bodies,
the Bethe ansatz goes actually back to Hermite, who solved
in this way the generalized Lam\'e equation, cf.\ \cite{WW}.
These integrable Schr\"odinger operators admit a $q$-deformation
due to Rujsenaars \cite{R}. 
In Sect. \ref{sixths} we present a $q$-deformation of
Hermite's result, i.e., we give eigenfunctions for the 
two-body Rujsenaars operator. 
These eigenfunctions are parametrized by a {\it spectral
curve}, similarly to the differential case. It is
a double covering of a hyperelliptic curve.
This result follows from our result and 
the observation that the transfer matrix associated to
an evaluation representation is equal to the Ruijsenaars
operator up to a scalar factor depending on the spectral
parameter. A similar observation, relating 
transfer matrices to Ruijsenaars operators,
has been made recently by Hasegawa \cite{H} in the context
of Sklyanin's algebra, indicating that our 
construction should be extendable to the $N$-body case.

\section{Modules over $E_{\tau,\eta}(sl_2)$ and transfer matrices}
\label{firsts}
In this paper, we construct eigenvectors of the transfer matrix 
of the ellipic quantum group $E_{\tau,\eta}(sl_2)$,
\cite{Fe,FV}, associated to certain highest weight modules.

Let us recall the definitions: we fix two complex parameters
$\tau,\eta$, such that Im$(\tau)>0$. 
The definition of $E_{\tau,\eta}(sl_2)$ is based on an $R$-matrix
$R(z,\lambda)$ which we now introduce. Let
\begin{equation}\label{jacobi}
 \theta(t)=-\sum_{j\in\Z}
e^{\pi i(j+\half)^2\tau+2\pi i(j+\half)(t+\half)},
\end{equation}
be Jacobi's theta function and
\be
\alpha(z,\lambda)=
\frac{\theta(\lambda+2\eta)\theta(z)}
{\theta(\lambda)\theta(z-2\eta)},\qquad
\beta(z,\lambda)=
-\frac{\theta(\lambda+z)\theta(2\eta)}
{\theta(\lambda)\theta(z-2\eta)},
\ee
Let $V$ be a two dimensional
complex vector space with basis $e[1],e[{-1}]$, and 
let $E_{ij}e[k]=\delta_{jk}e[i]$, $h=E_{11}-E_{-1,-1}$.
 Then, for $z,\lambda\in\C$, $R(z,\lambda)\in\End(V\otimes V)$ 
is the matrix
\bea\label{exR}
 R(z,\lambda)&\!=\!&E_{11}\otimes E_{11}+E_{-1,-1}\otimes E_{-1,-1}+
 \alpha(z,\lambda)E_{11}\otimes E_{-1,-1}
\\
 &
+\!&\alpha(z,-\lambda)E_{-1,-1}\otimes E_{11}
+\beta(z,\lambda)E_{1,-1}\otimes E_{-1,1}
+\beta(z,-\lambda)E_{-1,1}\otimes E_{1,-1}.
\eea
This $R$-matrix obeys the dynamical (or modified) quantum Yang--Baxter equation
\bea
R^{(12)}(z\!-\!w,\lambda\!-\!2\eta h^{(3)})\!\!\!\!\!
&R^{(13)}(z,\lambda)\,
R^{(23)}(w,\lambda\!-\!2\eta h^{(1)})&\!\!\!\!\!
=
\\
  &R^{(23)}(w,\lambda)\,
R^{(13)}(z,\lambda\!-\!2\eta h^{(2)})&\!\!\!\!\!
R^{(12)}(z\!-\!w,\lambda)
\eea
in $\End(V\otimes V\otimes V)$, $z,w,\lambda\in\C$.  The meaning of this
notation is the following: $R^{(12)}(\lambda-2\eta h^{(3)})v_1\otimes
v_2\otimes v_3$
is defined as
\be
(R(z,\lambda-2\eta\mu_3)v_1\otimes v_2)\otimes v_3,
\ee
if $h v_3=\mu_3v_3$.
The other terms are defined similarly: in general, let 
$V_1,\dots, V_n$ be modules over the one dimensional Lie
algebra $\h=\C h$ with one generator $h$, such that, for all $i$,
$V_i$ is the direct sum of finite dimensional eigenspaces
$V_i[\mu]$  of $h$, labeled by the eigenvalue $\mu$.
We call such modules diagonalizable $\h$-modules.
If $X\in\End(V_i)$
we denote by $X^{(i)}\in\End(V_1\otimes\dots\otimes V_n)$ 
the operator $\cdots\otimes\Id\otimes X\otimes\Id\otimes\cdots$
acting non-trivially on the $i$th factor, and
if $X=\sum X_k\otimes Y_k\in\End(V_i\otimes V_j)$ we set
$X^{(ij)}=\sum X_k^{(i)}Y_k^{(j)}$. If $X(\mu_1,\dots,\mu_n)$
is a function with values in $\End(V_1\otimes\dots\otimes V_n)$,
then $X(h^{(1)},\dots,h^{(n)})v=X(\mu_1,\dots,\mu_n)v$ if
$h^{(i)}v=\mu_iv$, for all $i=1,\dots,n$.

 Now, by definition, a module over $E_{\tau,\eta}(sl_2)$ is
a diagonalizable $\h$-module $W=\oplus_{\mu\in\C}  W[\mu]$,
together with an $L$-operator $L(z,\lambda)\in\End_\h(V\otimes W))$
(a linear map commuting with $h^{(1)}+h^{(2)}$)
depending meromorphically on  $z,\lambda\in\C$ and obeying the
relations
\bea
R^{(12)}(z\!-\!w,\lambda\!-\!2\eta h^{(3)})\!\!\!\!\!
&L^{(13)}(z,\lambda)\,
L^{(23)}(w,\lambda-2\eta h^{(1)})&\!\!\!\!\!
=
\\
  &L^{(23)}(w,\lambda)\,
L^{(13)}(z,\lambda\!-\!2\eta h^{(2)})&\!\!\!\!\!
R^{(12)}(z\!-\!w,\lambda)
\eea
For example, $W=V$, $L(w,\lambda)=R(w-z_0,\lambda)$ 
is a module over $E_{\tau,\eta}(sl_2)$,
called the fundamental representation, with evaluation point $z_0$.
In \cite{FV} more general examples of such modules were constructed: in 
particular, for any pair of complex numbers $\Lambda, z$ we
have an {\em evaluation Verma module} $V_\Lambda(z)$. 
Also, we have a notion of tensor products of modules over
$E_{\tau,\eta}(sl_2)$. The main examples considered in this
paper will be tensor products $V_{\Lambda_1}(z_1)\otimes
\cdots \otimes V_{\Lambda_n}(z_n)$ of evaluation Verma modules and
some of their subquotients.

For any module $W$ over $E_{\tau,\eta}(sl_2)$,  we define the
associated operator algebra, an algebra of operators on the space $\Fun(W)$
of meromorphic functions of $\lambda\in\C$
 with values in $W$. It is  generated by $h$,
acting on the values, and
operators  $a(z), b(z), c(z), d(z)$. Namely, let $\tilde L(z)
\in \End(V\otimes\Fun(W))$ be the operator
defined by 
$(\tilde L(z)(v\otimes f))(\lambda)=
L(z,\lambda)(v\otimes f(\lambda-2\eta\mu))$ 
if $h v=\mu v$. View $\tilde L(z)$
as a 2 by 2 matrix  with entries in $\End(\Fun(W))$:
\bea
\tilde L(z) (e[1]\otimes f)&\!=\!&e[1]\otimes a(z)f+e[{-1}]\otimes c(z)f,
\\
\tilde L(z) (e[{-1}]\otimes f)&\!=\!&e[1]\otimes b(z)f+e[{-1}]\otimes d(z)f.
\eea
The relations obeyed by these operators are described in detail
in \cite{FV} (in \cite{FV} these operators are denoted by 
$\tilde a(z), \tilde b(z)$ and so on).
\begin{thm}\label{one}\cite{FV}
For any module $W$ over $E_{\tau,\eta}(sl_2)$, the transfer
matrices $T(z)=a(z)+d(z)$ preserve the space
$H=\Fun(W)[0]$ of functions with values in the zero weight
space $W[0]$, and commute pairwise on $H$: $T(z)T(w)=T(w)T(z)$
on $H$.
\end{thm}

\begin{thm}\label{two}\cite{FV}
Let $W=V_{\Lambda_1}(z_1)\otimes\cdots\otimes V_{\Lambda_n}(z_n)$ be
a tensor product of evaluation modules, and let
 $\Lambda=\Lambda_1+\cdots+\Lambda_n$. Then $W[\Lambda]=\C v_0$,
where $v_0$, viewed as a constant function in $\Fun(W)$ obeys
the following highest weight condition: for every $z$,
 $c(z)v_0=0, a(z)v_0=A(z,\lambda)v_0,
d(z)v_0=D(z,\lambda)v_0$, with highest weight functions
\begin{equation}\label{AD}
A(z,\lambda)=1,\qquad
D(z,\lambda)=\frac{\theta(\lambda-2\eta\Lambda)}
{\theta(\lambda)}\prod_{j=1}^n
\frac{\theta(z-p_j)}
{\theta(z-q_j)},
\end{equation}
where we set $p_j=z_j+\eta(-\Lambda_j+1)$, $q_j=z_j+\eta(\Lambda_j+1)$
\end{thm}
A {\em highest weight module} $W$ of highest weight $(\Lambda,A,D)$,
is a module with a highest weight vector $v_0\in\Fun(W)$ such that
 $c(z)v_0=0, a(z)v_0=A(z,\lambda)v_0,
d(z)v_0=D(z,\lambda)v_0$, for all $z,\lambda$, and so that $\Fun(W)$
is spanned by the vectors of the form $b(t_1)\cdots b(t_j)v_0$,
as a vector space over the field of meromorphic functions of $\lambda$.
It is shown in \cite{FV} that if $A,D$ are of the form \Ref{AD},
for some $p_k,q_k$, then every irreducible highest weight module of
weight $(\Lambda,A,D)$ is isomorphic to a subquotient
of $\otimes V_{\Lambda_i}(z_i)$, where the parameters
$z_j,\Lambda_j$ are related to $p_k,q_k$ as in the theorem.
  If $p_k$, $q_k$ are generic under the condition that
$\sum (p_i-q_i)=2\eta\Lambda$, then
all highest weight modules of
weight $(\Lambda,A,D)$ are isomorphic to
$\otimes V_{\Lambda_i}(z_i)$ itself.

\section{Bethe ansatz}\label{seconds}

In this section we fix a highest weight module $W$ of weight
$(\Lambda,A,D)$ of the form \Ref{AD}, with highest weight vector
$v_0$.
We assume that $\Lambda$  is an even integer
$2m\geq 0$, so that the zero-weight space $W[0]$ can be nontrivial.
  
We follow the strategy of \cite{TF}: we  seek common eigenvectors of $T(w)$ 
in the form $b(t_1)\cdots b(t_m)v$, where $v\in \Fun(W)[\Lambda]$.
 The problem is to find conditions for
$t_1,\dots,t_n, v$ so that we have an eigenvector. The question
of completeness, i.e., whether ``all'' eigenvectors can be obtained
in this way, will not be addressed here, except in the example
of the $q$-analogue of the Lam\'e equation discussed below.

 Any non-zero vector $v\in\Fun(W)[\Lambda]$ is of the form 
$v=g(\lambda)v_0$, for some meromorphic function $g\neq 0$.
To find for which values of $t_1,\dots,t_m$ and which choice of $v$
we get an eigenvector
we commute $a$ and $d$ with the $b$'s using the relations of
the quantum group.
The relevant commutation relations are (see \cite{FV})
\bea
a(w)b(t)&\!=\!& r(t-w,\lambda)
b(t)a(w)
+ s(t-w,\lambda)
b(w)a(t),
\\
d(w)b(t)&\!=\!& r(w-t,\lambda-2\eta h)
b(t)d(w)
-s(t-w,\lambda-2\eta h)
b(w)d(t).
\eea
The functions $r,s$ are defined by the formulae
\be
r(t,\lambda)=\frac{\theta(t-2\eta)\theta(\lambda)}
     {\theta(t)\theta(\lambda-2\eta)},
\qquad
s(t,\lambda)=
\frac{\theta(t+\lambda)\theta(2\eta)}
     {\theta(t)\theta(\lambda-2\eta)}.
\ee
Note that these coefficients have poles for $t\in\Z+\tau\Z$, so
we have to be careful to avoid infinities in the calculations: we
assume that $w,t_1,\dots,t_m$ are all distinct modulo $\Z+\tau\Z$.
Using repeatedly the commutation relations to
bring $a$ to the right of the $b$'s, we get 
\bea
a(w)b(t_1)\cdots b(t_m)&\!=\!&
A_0
b(t_1)\cdots b(t_m)a(w)\\
 & &+
\sum_{j=1}^n A_jb(t_1)\cdots
b(t_{j-1})b(w)b(t_{j+1})\cdots b(t_m)a(t_j)
\eea
for some complex coefficients $A_j=A_j(w,t_1,\dots,t_m,\lambda)$.
The first term is called the ``wanted'' term, and the others
are ``unwanted'' terms. Similarly, we have
\bea
d(w)b(t_1)\cdots b(t_m)&\!=\!&
D_0
b(t_1)\cdots b(t_m)d(w)
\\
 & &+
\sum_{j=1}^n D_jb(t_1)\cdots
b(t_{j-1})b(w)b(t_{j+1})\cdots b(t_m)d(t_j)
\eea
for some complex coefficients $D_j=D_j(w,t_1,\dots,t_m,\lambda)$.
The coefficients $A_0$ and $A_1$ are easy to compute, since
they are  products of coefficients appearing in the commutation
relations:
\bea
A_0&\!=\!&\prod_{j=1}^mr(t_j-w,\lambda+2\eta(j-1)),
\\
A_1&\!=\!&s(t_1-w,\lambda)\prod_{j=2}^mr(t_j-t_1,\lambda+2\eta(j-1)).
\eea
A similar calculation gives
\bea
D_0&\!=\!&\prod_{j=1}^mr(w-t_j,\lambda-2\eta(j-1)),
\\
D_1&\!=\!&-s(t_1-w,\lambda)\prod_{j=2}^mr(t_1-t_j,\lambda-2\eta(j-1)).
\eea
A direct calculation of the other coefficients is more complicated.
However the answer is simple, thanks to the
\begin{lemma}\label{sym}
 For any permutation $\sigma$ of $m$ letters, and
all $j=1,\dots, m$,
\bea
A_j(w,t_{\sigma(1)},\dots,t_{\sigma(m)},\lambda)&\!=\!&
A_{\sigma(j)}(w,t_1,\dots,t_m,\lambda),\\
D_j(w,t_{\sigma(1)},\dots,t_{\sigma(m)},\lambda)&\!=\!&
D_{\sigma(j)}(w,t_1,\dots,t_m,\lambda),
\eea
\end{lemma}
The proof of this Lemma is deferred to Sect.\ \ref{sproofs}.

The next step is to find conditions for cancellation of unwanted
terms.
We have, with $v(\lambda)=g(\lambda)v_0$,
\bean\label{C}
T(w)b(t_1)\cdots b(t_m)v&\!=\!&
C_0
b(t_1)\cdots b(t_m)v
\\
 & &+
\sum_{j=1}^n C_jb(t_1)\cdots
b(t_{j-1})b(w)b(t_{j+1})\cdots b(t_m)v\notag
\eean
for some coefficients $C_j$. The condition of cancellation is
$C_j=0$, $j\geq 1$. Let us first consider $C_1$. It is given by
\bea
C_1(w,t,\lambda)&\!=\!&A_1(w,t,\lambda)
\frac{g(\lambda+2\eta(m-1))}{g(\lambda+2\eta m)}\\
 & &
+D_1(w,t,\lambda)D(t_1,\lambda+2\eta m)
\frac{g(\lambda+2\eta(m+1))}{g(\lambda+2\eta m)}.
\eea
The condition $C_1=0$ is then equivalent to
\bea
 &r(t_2\!-\!t_1,\lambda+2\eta)\cdots
r(t_m\!-\!t_1,\lambda+2\eta(m-1))g(\lambda+2\eta(m-1))&
\\
=
 &r(t_1\!-\!t_2,\lambda-2\eta)\cdots
r(t_1\!-\!t_m,\lambda-2\eta(m-1))D(t_1,\lambda+2\eta m)g(\lambda+2\eta(m+1))&
\eea
which can be written as
\be
\prod_{j=2}^m
\frac
{\theta(t_j-t_1-2\eta)}
{\theta(t_j-t_1+2\eta)}
\prod_{k=1}^n
\frac
{\theta(t_1-q_k)}
{\theta(t_1-p_k)}
\frac
{\theta(\lambda+2\eta(m-1))\theta(\lambda+2\eta m)
g(\lambda+2\eta(m-1))}
{\theta(\lambda)
\theta(\lambda-2\eta)g(\lambda+2\eta(m+1))}
=1
\ee
In particular, the left-hand side of this equation should be independent
of $\lambda$. This holds if $g$ is taken in the form
\be
g(\lambda)=e^{c\lambda}\prod_{j=1}^m
\frac{\theta(\lambda-2\eta j)}{\theta(2\eta)}.
\ee
Thus $C_1$ vanishes if $g$ has this form and the $t_i$ obey
the equation
\be
\prod_{j=2}^n\frac
{\theta(t_j-t_1-2\eta)}
{\theta(t_j-t_1+2\eta)}
\prod_{k=1}^n
\frac
{\theta(t_1-q_k)}
{\theta(t_1-p_k)}
=
e^{4\eta c}
\ee
Using Lemma \ref{sym} we find the conditions for $C_j$ to 
vanish for $j=1,\dots, m$. The result is:

\begin{thm}\label{three}
Let $W$ be a highest weight module over $E_{\tau,\eta}(sl_2)$
of highest weight $(\Lambda,A,D)$ with $A, D$ of the form
\Ref{AD} and $\Lambda=2m\in2\Z_{\geq 0}$. Let $T(w)\in\End(H)$, $w\in\C$,
 be
the corresponding transfer matrices. Let
 $v(\lambda)=e^{c\lambda}\prod_{j=1}^m\theta(\lambda-2\eta j)v_0$.
Then, for any solution $(t_1,\dots,t_m)$ of the Bethe ansatz
equations
\begin{equation}\label{bae}
\prod_{j:j\neq i}
\frac
{\theta(t_j-t_i-2\eta)}
{\theta(t_j-t_i+2\eta)}
\prod_{k=1}^n
\frac
{\theta(t_i-q_k)}
{\theta(t_i-p_k)}
=
e^{4\eta c},\qquad i=1,\dots,m,
\end{equation}
such that, for all $i<j$, $t_i\neq t_j\mod \Z+\tau\Z$,
the vector $b(t_1)\cdots b(t_m)v$ is a common eigenvector
of all transfer matrices $T(w)$ with eigenvalues
\be
\epsilon(w)=
e^{-2\eta c}\prod_{j=1}^m
\frac
{\theta(t_j-w-2\eta)}
{\theta(t_j-w)}
+
e^{2\eta c}\prod_{j=1}^m
\frac
{\theta(t_j-w+2\eta)}
{\theta(t_j-w)}
\prod_{k=1}^n
\frac
{\theta(w-p_k)}
{\theta(w-q_k)}
\ee
\end{thm}
What is left to prove is the formula for the eigenvalue, which is
given by $C_0$ in \Ref{C}. By definition 
\be
C_0=A_0
\frac
{g(\lambda+2\eta(m-1))}
{g(\lambda+2\eta m)}
+D_0
\frac
{g(\lambda+2\eta(m+1))}
{g(\lambda+2\eta m)}
D(w,\lambda+2\eta m).
\ee
By inserting the formulas for $A_0,D_0,g,D$, we see that the 
$\lambda$ dependence disappears for the same reason as before,
and we are left with the formula for $\epsilon$ given in the
theorem.

We conclude this section by giving an explicit formula
for $b(t_1)\cdots b(t_m)v$. It is sufficient to consider
the case where $W=V_{\Lambda_1}(z_1)\otimes\cdots\otimes V_{\Lambda_n}
(z_n)$, for any other highest weight module considered in
the previous Theorem is a quotient of the submodule 
of a module of this form generated by the product of highest
weight vectors. Recall \cite{FV} that $V_{\Lambda}(z)$ is defined
to be the infinite dimensional vector space with basis $e_0,e_1,e_2,\dots$,
such that $h e_j=(\Lambda-2j)e_j$, and with the action of the
generators of $E_{\tau,\eta}(sl_2)$ given by explicit formulae.

\begin{thm}\label{four}
With the notations as in the previous theorem,
\bea
b(t_1)\cdots b(t_m)v&\!\!\!=\!\!\!&
(-1)^m
e^{c(\lambda+2\eta m)}
\sum_{I_1,\dots,I_n}
\prod_{l=1}^n
\prod_{i\in I_l}
\prod_{k=l+1}^n
\frac
{\theta(t_i-p_k)}
{\theta(t_i-q_k)}
\\
 & &\!\!\!\!\!\!\times\prod_{k<l}
\prod_{i\in I_k,j\in I_l}
\frac
{\theta(t_i-t_j-2\eta)}
{\theta(t_i-t_j)}\\
 & &\!\!\!\!\!\!\times
\prod_{k=1}^{n}
\prod_{j\in I_k}
\frac
{\theta(\lambda\!+\!t_j\!-\!q_k\!+\!2\eta
m_k\!-\!2\eta\sum_{l=k+1}^n(\Lambda_l\!-\!2m_l))}
{\theta(t_j-q_k)}
e_{m_1}\!\otimes\cdots\otimes\! e_{m_n}.
\eea
The summation is over all partitions of
$\{1,\dots,m\}$ into $n$ disjoint subsets $I_1,\dots,I_n$. 
The cardinality of $I_j$
is denoted by $m_j$.
\end{thm}
In  particular, we see that our eigenvector $\psi(\lambda)$
is an entire function 
of $\lambda$ obeying $\psi(\lambda+1)=(-1)^me^c\psi(\lambda)$.
The proof of this Theorem is contained in Section \ref{sproofs}.

\section{Discrete models}\label{thirds}

The construction of the transfer matrix admits the following
variation. The difference operators $a(w),\dots,d(w)$
shift the argument of functions by $\pm 2\eta$. Therefore we may replace
$\Fun(W)$ by the space $\Fun_\mu(W)$ of all functions from the set
$C_\mu=\{\mu+2\eta j,j\in\Z\}$ to $W$. The operators are then well-defined
on $\Fun_\mu(W)$ if $\mu$ is generic. Also, it follows from
Theorem \ref{four} that the restriction to $C_\mu$ of the
Bethe ansatz eigenfunctions is well defined for all $\mu$.
We thus have:

\begin{corollary}
Suppose $t_1,\dots,t_m$ is a solution to the Bethe ansatz equations
\Ref{bae}. Then,
for generic $\mu$, the restriction to $C_\mu$ of $b(t_1)\cdots b(t_m)v$
is a common eigenfunction of the operators
$T(w)\in\End(\Fun_\mu(W)[0])$.
\end{corollary}

\section{Interaction-round-a-face models}\label{fourths}

In this section, we consider a special case of the above construction,
and relate $T(w)$ to the transfer matrix of the interaction-round-a-face
(IRF) (also called solid-on-solid) 
models of \cite{Ba}, \cite{ABF}. Therefore, our formulae give in particular
eigenvectors for transfer matrices of IRF models.  

Notice first (see \cite{Fe}) that if we define%
\footnote{We adopt here a convention for the definition of the 
Boltzmann weights which is slightly different than in \cite{Fe} and in
better agreement with the literature.  }
a ``Boltzmann weight''
$w(a,b,c,d;z)$, depending on complex parameters $a,b,c,d,z$, such that
$a-b,b-c,c-d,a-d\in\{1,-1\}$, by \be R(z,-2\eta d)e[{c-d}]\otimes e[{b-c}]
=\sum_a w(a,b,c,d;z)e[{b-a}]\otimes e[{a-d}], \ee (the sum is over one
or two allowed values of $a$) then the dynamical quantum Yang--Baxter
equation translates into the star-triangle relation \bea &\sum_g
w(a,b,g,f;z-w) 
w(f,g,d,e;z) 
w(g,b,c,d;w)
& \\ &= 
\sum_g 
w(f,a,g,e;w)
w(a,b,c,g;z) 
w(g,c,d,e;z-w).& 
\eea 
We let $W=V^{\otimes n}$ be the
tensor product of fundamental representations with evaluation points
$z_1,\dots,z_n$. Then \be
L(z,\lambda)=R^{(01)}(z-z_1,\lambda-2\eta\sum_{j=2}^{n}h^{(j)}) \cdots
R^{(0,n-1)}(z-z_{n-1},\lambda-2\eta h^{(n)}) R^{(0n)}(z-z_n,\lambda).
\ee In this formula the factors of $V$ in $V\otimes W=V^{\otimes{n+1}}$
are numbered from $0$ to $n$. The module $W$ is a highest weight module
with highest weight functions of the form \Ref{AD} with $p_j=z_j$,
$q_j=z_j+2\eta$.

\begin{figure}
\begin{picture}(350,100)
\put(10,48){\makebox{$z$}}
\put(20,50){\line(1,0){310}}
\put(30,60){\makebox{$b_1$}}
\put(30,35){\makebox{$a_1$}}
\put(47,20){\makebox{$z_n$}}
\put(50,30){\line(0,1){40}}
\put(60,60){\makebox{$b_n$}}
\put(60,35){\makebox{$a_n$}}
\put(77,20){\makebox{$z_{n-1}$}}
\put(80,30){\line(0,1){40}}
\put(160,60){\makebox{$\cdots$}}
\put(160,35){\makebox{$\cdots$}}
\put(237,20){\makebox{$z_3$}}
\put(240,30){\line(0,1){40}}
\put(250,60){\makebox{$b_3$}}
\put(250,35){\makebox{$a_3$}}
\put(267,20){\makebox{$z_2$}}
\put(270,30){\line(0,1){40}}
\put(280,60){\makebox{$b_2$}}
\put(280,35){\makebox{$a_2$}}
\put(297,20){\makebox{$z_1$}}
\put(300,30){\line(0,1){40}}
\put(310,60){\makebox{$b_1$}}
\put(310,35){\makebox{$a_1$}}
\end{picture}
\caption{Graphical representation of the row-to-row transfer matrix
of an IRF model. Each crossing represents a Boltzmann weight $w$
whose arguments are the labels of the adjoining regions and
the difference of the parameters associated to the lines.}
\label{fig1}
\end{figure}
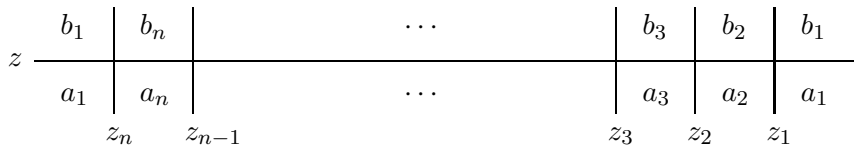

Let us now introduce a basis $|a_1,\dots,a_n\rangle$ of
$\Fun_\mu(W[0])$, 
labeled by $a_i\in \mu+\Z$ with $a_i-a_{i+1}\in\{1,-1\}$,
$i=1,\dots,n-1$, and $a_n-a_1\in\{1,-1\}$. We let $\delta(\lambda)=
1$ if $\lambda=0$ and $0$ otherwise. Then we define
\be
|a_1,\dots,a_n\rangle(\lambda)
=\delta(\lambda+2\eta a_1)e[{a_1-a_2}]\otimes e[{a_2-a_3}]\otimes
\cdots\otimes e[{a_n-a_1}].
\ee
If $\Gamma$ is the shift operator $\Gamma f(\lambda)=f(\lambda-2\eta)$,
then $\Gamma|a_1,\dots,a_n\rangle=|a_1\!-\!1,\dots,a_n\!-\!1\rangle$.
Using this, and the fact that $h^{(j)}|a_1,\dots,a_n\rangle=
(a_{j+1}\!-\!a_{j})|a_1,\dots,a_n\rangle$, we get
\be
T(z)|a_1,\dots,a_n\rangle
=\sum_{b_1,\dots,b_n}\prod_{j=1}^n
w(b_j,a_j,a_{j+1},b_{j+1};z\!-\!z_j)|b_1,\dots,b_n\rangle,
\ee
with the understanding that $b_{n+1}=b_1$, $a_{n+1}=a_1$. The
(finite) sum is over the values of the indices $b_i$ for which
the Boltzmann weights are defined.
Comparing with \cite{Ba}, we see that $T(z)$,
in this basis, is the row-to-row transfer matrix of the
(inhomogeneous) interaction-round-a-face model
associated to the solution $w(a,b,c,d;z)$ of the star-triangle equation
(see \cite{Ba}). The situation is best visualized by 
looking at the graphical representation of Fig.\ \ref{fig1}.

This construction can be in principle extended to higher
representations, and we obtain in this way eigenvectors of
transfer matrices of the IRF models of \cite{DJMO}.

\section{The eight-vertex model}\label{fifths}

We show in this section how to obtain from our result 
eigenvectors for the
transfer matrix of the eight-vertex model. 

The eight-vertex model is based on Baxter's solution
\bea
R_{8V}(z)=
a_{8V}(z)(E_{1,1}\otimes E_{1,1}+E_{-1,-1}\otimes E_{-1,-1})\\
+
b_{8V}(z)(E_{1,-1}\otimes E_{1,-1}+E_{-1,1}\otimes E_{-1,1})\\
+
c_{8V}(z)(E_{1,-1}\otimes E_{-1,1}+E_{-1,1}\otimes E_{1,-1})\\
+
d_{8V}(z)(E_{1,1}\otimes E_{-1,-1}+E_{-1,-1}\otimes E_{1,1}),
\eea
of the Yang--Baxter equation. The coefficients are
\bea
a_{8V}(z)=\frac{\theta_0(z)\theta_0(2\eta)}{\theta_0(z-2\eta)\theta_0(0)},
\qquad
b_{8V}(z)=\frac{\theta_1(z)\theta_0(2\eta)}{\theta_1(z-2\eta)\theta_0(0)},
\\
c_{8V}(z)=-\frac{\theta_0(z)\theta_1(2\eta)}{\theta_1(z-2\eta)\theta_0(0)},
\qquad
d_{8V}(z)=-\frac{\theta_1(z)\theta_1(2\eta)}{\theta_0(z-2\eta)\theta_0(0)},
\eea
in terms of the theta functions with characteristics
\be
\theta_1(z)=-\sum_{j\in\Z}e^{2\pi i(j+\half)^2\tau+2\pi i(j+\half)(z+\half)},
\qquad
\theta_0(z)=ie^{-i\pi i(z+\tau/2)}\theta_1(z).
\ee
{\it Remarks.} Baxter uses $-z$ instead of $z$. The classical Jacobi notation,
used by Baxter, is $\theta_1(z)=\mbox{H}(2Kz)$, $\theta_0(z)=\Theta(2Kz)$.

\medskip
For any (generic) $z_1,\dots,z_n$, one then defines commuting
 transfer matrices 
\begin{equation}\label{t8v}
T_{8V}(z)=\tr_0R_{8V}(z-z_1)^{(01)}\cdots R_{8V}(z-z_n)^{(0n)}.
\end{equation}
acting on $(\C^2)^{\otimes n}$.
 The
relation with our transfer matrices is based on the following
identity, which is a version of Baxter's vertex-IRF transformation:

\begin{proposition}\label{l7}
Let $S(z,\lambda)$ be the matrix
\be
\frac1{\theta(\lambda)}\left(
\begin{matrix}
\theta_0(z-\lambda+1/2) & -\theta_0(-z-\lambda+1/2)\\
-\theta_1(z-\lambda+1/2) & \theta_1(-z-\lambda+1/2)
\end{matrix}\right).
\ee
Then 
\be
S(w,\lambda)^{(2)}S(z,\lambda-2\eta h^{(2)})^{(1)}
R(z-w,\lambda)
=
R_{8V}(z-w)S(z,\lambda)^{(1)}S(w,\lambda-2\eta h^{(1)})^{(2)}
\ee
\end{proposition}
\begin{prf} We need to recall some well-known facts about
the functions $\theta_\alpha$ and Baxter's $R$-matrix $R_{8V}$.

The functions $\theta_\alpha$, $\alpha=0,1$ are entire, and uniquely
determined up to normalization 
by the properties $\theta_\alpha(z+1)=(-1)^\alpha\theta_\alpha(z)$,
$\theta_\alpha(z+\tau)=ie^{-\pi i(z+\tau/2)}\theta_{1-\alpha}(z)$.
 Also, $\theta_\alpha(-z)=(-1)^\alpha\theta_\alpha(z)$,
and the zeros of $\theta_\alpha$ are simple and of the form
$r+2s\tau$ if $\alpha=1$ and $r+(2s+1)\tau$ if $\alpha=0$, ($r,s\in\Z$).
We have $\theta_0(z)\theta_1(z)=C(\tau)\theta(z)$ for some
constant $C(\tau)$, as one can see comparing transformation properties
under translations by $\Z+\tau\Z$.

Let $A=\left(\begin{matrix}-1&0\\0&1\end{matrix}\right)$,
$B=\left(\begin{matrix}0&1\\1&0\end{matrix}\right)$.
Then $R_{8V}(z)$ commutes with $A\otimes A$ and $B\otimes B$.
Moreover $R_{8V}(z+1)=A^{(1)}R_{8V}(z)A^{(1)}$, 
$R_{8V}(z+\tau)=e^{-2\pi i\eta}B^{(1)}R_{8V}(z)B^{(1)}$, and $R_{8V}(0)=P$,
the flip $u\otimes v\mapsto v\otimes u$. As a function of $z$,
$R_{8V}$ is meromorphic. Its poles are simple and are at $2\eta$
modulo $\Z+\tau\Z$. The residue of $R_{8V}(z)$
at $z=2\eta$ is $\theta(2\eta)/\theta'(0)$ times the anstisymmetrization
operator $\Pi: u\otimes v\mapsto u\otimes v-v\otimes u$.

Let $\phi(z)=(\theta_1(z),\theta_0(z))$. This vector is the
(up to normalization) unique vector with entire holomorphic components
such that $\phi(z+1)=A\phi(z)$, $\phi(z+\tau)=ie^{-i\pi(z+\tau/2)}
B\phi(z)$.

We claim that the statement of the Proposition is equivalent to
the following set of identities, which, incidentally, is 
essentially the
standard form of the vertex-IRF transformation.

\begin{lemma}
Let $\phi^\pm(z,\lambda)=\phi(\mp z-\lambda+\half)$. Then
\bea 
R_{8V}(z-w)
\phi^{\pm}(z,\lambda\mp 2\eta)\otimes
\phi^{\pm}(w,\lambda)
&\!=\!&
\phi^{\pm}(z,\lambda)\otimes
\phi^{\pm}(w,\lambda\mp 2\eta)\\
R_{8V}(z-w)
\phi^{\pm}(z,\lambda\pm 2\eta)\otimes
\phi^{\mp}(w,\lambda)
&\!=\!&
\alpha(z-w,\pm\lambda)
\phi^{\pm}(z,\lambda)\otimes
\phi^{\mp}(w,\lambda\mp 2\eta)
\\
 & &+
\beta(z-w,\pm\lambda)
\phi^{\mp}(z,\lambda)\otimes
\phi^{\pm}(w,\lambda\pm 2\eta).
\eea
\end{lemma}
The (known) proof of this lemma consists in comparing transformation
properties under lattice translations and poles of both sides
of the equations as functions of $z$, $w$, and using the uniqueness
of $\phi$.

It remains to show that these identities are equivalent to
the statement of the Proposition. Let us write them in matrix form: 
Let $\hat S(z,\lambda)$ be the 2 by 2 matrix whose columns are
$\phi^+$ and $\phi^-$ then the previous lemma reads
\be
R_{8V}(z-w)\hat S(w,\lambda)^{(2)}\hat S(z,\lambda-2\eta h^{(2)})^{(1)}
=
\hat S(z,\lambda)^{(1)}\hat S(w,\lambda-2\eta h^{(1)})^{(2)}
R(z-w,\lambda)^t.
\ee
The transposed of a matrix in $\End(\C^2\otimes C^2)$ is defined
by the rule $(X\otimes Y)^t=X^t\otimes Y^t$, $X,Y\in\End(\C^2)$.
The statement of the Proposition is proved by transposing both
sides of the equation, and inverting the matrices $\hat S$. The
identity we need is $(\hat S^t)^{-1}=\mbox{const}\,S$,
which follows from the identity $\det \hat S(z,\lambda)=
\mbox{const}\,\theta(z)\theta(\lambda)$. The latter formula
can be proved by comparing the transformation properties under
lattice translations of $z$ and $\lambda$, and using the fact that
$\theta$ is the unique function (up to normalization) such
that $\theta(z+1)=-\theta(z)$ and $\theta(z+\tau)=
-e^{-\pi i(\tau+2z)}\theta(z)$.
\end{prf}

The transfer matrix $T(z)$ of the previous section has the
form $T(z)\psi(\lambda)=a(z,\lambda)\psi(\lambda-2\eta)
+d(z,\lambda)\psi(\lambda+2\eta)$, where $a(z,\lambda)$,
 $d(z,\lambda)\in\End((\C^2)^{\otimes n}[0])$ are the diagonal
matrix elements of $L(z,\lambda)$.

\begin{corollary} Fix generic $z_1,\dots,z_n\in\C$.
Let 
\be S_n(\lambda)=S(z_n,\lambda)^{(n)}
S(z_{n-1},\lambda-2\eta h^{(n)})^{(n-1)}
\cdots
S(z_1,\lambda-2\eta\sum_{j\geq 2}h^{(j)})^{(1)}.\ee
 Then,
for all generic complex $z,\lambda$,
\be
T_{8V}(z)S_n(\lambda)=S_n(\lambda+2\eta)a(z,\lambda+2\eta)
+S_n(\lambda-2\eta)d(z,\lambda-2\eta),
\ee
on the zero weight space $(\C^2)^{\otimes n}[0]$.
\end{corollary}
\begin{prf}
Let us write $T_{8V}(z)=\tr_0L_{8V}(z)$, see \Ref{t8v}. The
matrix $L_{8V}(z)$ acts on the tensor product of $n+1$ copies
of $\C^2$, numbered from $0$ to $n$. By iterating Proposition \ref{l7},
we obtain
\begin{equation}\label{SLLS}
S_n(\lambda)
S(z,\lambda-2\eta h)^{(0)}L(z,\lambda)=
L_{8V}(z)S(z,\lambda)^{(0)}S_n(\lambda-2\eta h^{(0)}).
\end{equation}
In this formula $S_n(\lambda)$ acts on the factors numbered 
from 1 to $n$,
and $h=h^{(1)}+\cdots+h^{(n)}$. Let $L_\alpha^\beta(z,\lambda)$
be defined by 
\be
L(z,\lambda)e[\alpha]=\sum_{\beta=\pm1} e[\beta]\otimes
L_\alpha^\beta(z,\lambda),\ee
and set $\psi_\alpha=S(z,\lambda+2\eta\alpha)e[\alpha]$.
Replacing $\lambda$ by
$\lambda+2\eta\alpha$ in \Ref{SLLS}
and acting on a vector of the form  $e[\alpha]\otimes u$,
$\alpha\in\{1,-1\}$, where $hu=0$, yields
\be
\sum_\beta 
\psi_\beta\otimes S_n(\lambda+2\eta\alpha)L_\alpha^\beta(z,\lambda
+2\eta\alpha)u=
L_{8V}(z)(
\psi_\alpha\otimes S_n(\lambda)u),
\ee
where we used the fact that $L(z,\lambda)$ commutes with
$h^{(0)}+h$. Since, for generic parameters,
$\psi_1, \psi_{-1}$ form
a basis of $\C^2$, and $L_1^1(z,\lambda)=a(z,\lambda)$,
$L_{-1}^{-1}(z,\lambda)=d(z,\lambda)$, the proof is 
complete.
\end{prf}

\begin{thm}
Let $f\mapsto \int f(\lambda)$ be a linear function on a
space $\cal F$ of functions of $\lambda\in\C$, such that
$\int f(\lambda+2\eta)=\int f(\lambda)$ for all $f\in{\cal F}$.
Extend $\int$ to vector-valued functions by acting componentwise.
Then for each eigenfunction $\psi(\lambda)$ of $T(z)$, the vector
$\int S_n(\lambda)\psi(\lambda)$, if defined, is an eigenvector of
$T_{8V}(z)$ with the same eigenvalue.
\end{thm}

This theorem is an easy consequence of the previous corollary.
What is left to do is to find linear forms $\int$ defined
on the components of $S(\lambda)\psi(\lambda)$ for the Bethe
ansatz eigenfunctions $\psi$ of Theorems \ref{three}, \ref{four}.

This can be done easily in two situations: notice that both
$S(\lambda)$ and $\psi(\lambda)$ are periodic in $\lambda$
with period $2$, if the parameter $c$ belongs to $\pi i\Z$.
\noindent 1. If $\eta=p/q$ is rational, we may choose
\be
\int f(\lambda)=\sum_{j=0}^{q-1} f(\mu+2\eta j),
\ee
for any generic $\mu$.

\noindent 2. If $\eta$ is real, we set
\be
\int f(\lambda)=\int_0^2f(\mu+\lambda)d\lambda,
\ee
for generic $\mu$.

\section{$q$-deformed Lam\'e equation}\label{sixths}

We consider here the special case where $W$ is the 
evaluation module $V_{2m}(0)$
with positive integer $m$.
Then the zero weight space is one-dimensional.
{}From the expression for $a$ and $d$ given in
\cite{FV}, we
see that the transfer matrix has the form (cf.\ \cite{H}, where
a similar observation is made in the context of the Sklyanin algebra) 
\begin{equation}\label{eTL}
 T(z)=
\frac{\theta(z-\eta)}
{\theta(z-(2m+1)\eta)}L,
\end{equation}
where the difference operator $L$ is independent of $z$ and
is given by 
\be
L\psi(\lambda)=
\frac{\theta(\lambda+2\eta m)}
{\theta(\lambda)}
\psi(\lambda-2\eta)
+
\frac{\theta(\lambda-2\eta m)}
{\theta(\lambda)}
\psi(\lambda+2\eta).
\ee
This operator is in fact conjugated to the two-body Ruijsenaars operator
\cite{R}. 
It obeys $(L\phi,\psi)=(\phi,L\psi)$ with respect to
the symmetric bilinear form
\be
(\phi,\psi)=\int \frac{\phi(\lambda)\psi(-\lambda)}
{\prod_{j=1}^m\theta(\lambda-2\eta j)
\theta(\lambda+2\eta j)},
\ee
where $\int$ is a linear function invariant under shift by $2\eta$,
and change of sign of the argument,
 defined on a suitable space of functions. 
If $\eta$ is real, we may take $\int$ to be
the integral over 
$\gamma=\gamma_++\gamma_-$, where the straight path $\gamma_+$ 
goes from $\tau/2$
to $1+\tau/2$ and $\gamma_-$ goes from $-1-\tau/2$ to $-\tau/2$.
This linear form is defined, say, on the space of meromorphic
functions whose poles are contained in $\{2\eta j+k+l\tau,j,k,l\in\Z\}$,
which is preserved by the action of $L$.

Note that $L$ has periodic coefficients, and therefore preserves the
space of Bloch functions $\psi$ such that
$\psi(\lambda+1)=\mu\psi(\lambda)$.
The bilinear form is well-defined  on this space of functions.

We consider the eigenvalue problem 
\begin{equation}\label{lame}
L\psi=\epsilon\psi.
\end{equation}
It is
a $q$-deformation of the generalized Lam\'e equation: 
we have the expansion of $L$ in powers of $\eta$:
\be
L_\eta=2\,\Id+4\eta^2(\frac{d^2}{d\lambda^2}-2m\frac{\theta'(\lambda)}
{\theta(\lambda)}+m^2\frac{\theta''(\lambda)}{\theta(\lambda)})+O(\eta^4).
\ee
The differential operator appearing in the second order
coefficient is the generalized Lam\'e operator (up to conjugation
by $\theta(\lambda)^m$). Our Bethe ansatz solution
is a $q$-deformation of Hermite's solution of the (generalized)
Lam\'e equation (see the last pages of \cite{WW}).

\begin{thm}
Let $(t_1,\dots,t_m,c)$ be a solution of the Bethe ansatz
equations:
\begin{equation}\label{baelame}
\frac
{\theta(t_i-\eta(1+2m))}
{\theta(t_i-\eta(1-2m))}
\prod_{j:j\neq i}
\frac
{\theta(t_j-t_i-2\eta)}
{\theta(t_j-t_i+2\eta)}
=
e^{4\eta c},\qquad i=1,\dots,m,
\end{equation}
such that $t_i\neq t_j\mod \Z+\tau\Z$ if $i\neq j$. Then
\begin{equation}\label{psi3}
\psi(\lambda)=e^{c\lambda}\prod_{j=1}^m
{\theta(\lambda+t_j-\eta)},
\end{equation}
is a solution of the $q$-deformed Lam\'e equation $L\psi=\epsilon\psi$
with eigenvalue
\be
\epsilon=e^{-2\eta c}\frac{\theta(4\eta m)}{\theta(2\eta m)}
\prod_{j=1}^m
\frac{\theta(t_j+(2m-3)\eta)}
{\theta(t_j+(2m-1)\eta)}.
\ee
\end{thm}
\begin{prf} This theorem is a special case of Theorem \ref{three}.
It follows from \Ref{eTL} that the formula for the eigenvalue is $\epsilon=
\frac{\theta(z-(2m+1)\eta)}
{\theta(z-\eta)}
\epsilon(z)$, where $\epsilon(z)$ is the function in Theorem \ref{three}.
Since this expression is independent of $z$, we can evaluate it
at any $z$. The formula given in this theorem is obtained by taking
$z=(1-2m)\eta$, so that the second term in $\epsilon(z)$ vanishes.
\end{prf}

The Bethe ansatz equations have the form $b_i(t)=e^{4\eta c}$,
$i=1,\dots,m$. We may eliminate $c$ and consider the set of
$(t_1,\dots,t_m)$ obeying the equations $b_i(t)b_j(t)^{-1}=1$.
The functions $b_i(t)b_j(t)^{-1}$ are doubly periodic meromorphic
functions with periods 1 and $\tau$ in each of the variables $t_j$.
Also, eigenfunctions corresponding to solutions related by
shifts of the variables $t_i$ by 1 or $\tau$ are proportional.
Therefore the set
\be
X=\{t\in(\C/\Z+\tau\Z)^m/S_m| b_i(t)=b_j(t)\neq\infty, 
i=1,\dots,m,\; t_i\neq t_j
(i\neq j)\}
\ee
is an algebraic subvariety of the symmetric power of our elliptic
curve, which parametrizes our eigenfunctions. We call it the
Hermite--Bethe variety.
On this variety, we have a single-valued function $\epsilon^2$, the
square of the eigenvalue, which is the restriction to $X$ of a
doubly periodic function. Therefore it is an algebraic function
$\epsilon^2:X\to\C$. The eigenvalue, which is the square root of this function,
is two-valued on $X$. This is connected to the fact that generically
there are two eigenfunctions corresponding to a single point of
the Hermite--Bethe variety: a point of $X$ determines $c$ only
modulo $(4\eta)^{-1}2\pi i\Z$. If
$\psi$ is an eigenfunction corresponding to a point of $X$,
 then $e^{\pi i\lambda/2\eta}\psi(\lambda)$ is 
an eigenfunction with the opposite eigenvalue
corresponding to the same point but with
$c$ translated by $2\pi i/4\eta$. 

Note 
that the space of meromorphic solutions of the difference equation
\Ref{lame} is a vector space over the field $K$ of
meromorphic $2\eta$-periodic functions of $\lambda$. 
If $(t_1,\dots,t_m,c)$ is a solution of the Bethe ansatz 
equations, then, for all $k\in\Z$,  $(t_1,\dots,t_m,c+2\pi ik/2\eta)$
is also a solution with the same value of $\epsilon$. The
corresponding eigenfunction is proportional (with a $2\eta$-periodic
coefficient) to the original eigenfunction. Therefore
we should consider $c$ modulo $\pi i\eta^{-1}\Z$.

Note also that we have an involutive automorphism
$\sigma$ of $X$, sending $(t_1,\dots,t_m)$
to $(2\eta-t_1,\dots,2\eta-t_m)$. If $(t,c)$ is a solution
of the Bethe ansatz equations, then also $(\sigma(t),-c)$. The
corresponding eigenfunctions are related by the ``Weyl reflection''
$\psi(\lambda)\mapsto\psi(-\lambda)$ and have the same eigenvalue.

We now turn to the question of completeness.  Let $S^mE_\tau=
(\Z+\tau\Z)^m/S_m$ be the symmetric power of the elliptic curve.

\begin{thm} Suppose that $\eta\in\C$ is generic.
For generic $\epsilon\in\C$, there are precisely two  solutions
\be
(t_1,\dots,t_m,c)\qquad\mbox{and}\qquad
(2\eta-t_1,\dots,2\eta-t_m,-c), 
\ee
in $S^mE_\tau\times(\C/\pi i{\eta}^{-1}\Z)$ of the Bethe ansatz
equations \Ref{baelame} with given $\epsilon$. The corresponding
eigenfunctions $\psi_\pm$ are linearly independent over the field $K$ of
$2\eta$-periodic meromorphic functions of $\lambda$, and all
solutions of the
$q$-Lam\'e equation \Ref{lame} are linear combinations of 
$\psi_+$, $\psi_-$ with coefficients in $K$. 
\end{thm}

\begin{prf} It is straightforward to check that the reflection 
$t_i\mapsto 2\eta-t_i$, $c\mapsto -c$,
 maps solutions to solutions preserving
the value of the eigenvalue $\epsilon$. 

Let $\bar X$ be the closure of $X$ in the symmetric power of the 
elliptic curve. Then $\bar X$ contains the point 
$P=((1-2m)\eta,\dots, -3\eta, -\eta)$. In a neighborhood of
$P$ we can introduce local coordinates 
$u_1=t_1-(1-2m)\eta$, $u_j=t_j-t_{j-1}-2\eta$, $2\leq j\leq m$.
In these coordinates, $\bar X$ is described by the equations
$u_j=c_ju_m^{m-j+1}(1+O(u_m))$, $j=1,\dots,m-1$, for some constants
$c_j$. Therefore, in this neighborhood, $\bar X$ is a non-singular
curve, and $u_m$ is a local parameter at $P$. The eigenvalue has
the form 
\be
\epsilon=\mbox{const}\, u_m^{-1/2}(1+O(u_m))
\ee
in a neighborhood of $P$. Therefore $\epsilon$ is a non-constant
function on $X$. Similarly one shows that 
$c$ appearing in the Bethe ansatz equation is a non-constant
function on $X$ ($c$ diverges at $P$).
Since $\epsilon^2$ is algebraic, we have, for any generic value
of $\epsilon$, a solution $(t_1,\dots,t_m,c)$
of the Bethe ansatz equations, and thus an
eigenfunction $\psi_+$. Let $\psi_-$ be the eigenfunction with
the same eigenvalue associated
to the reflected solution $(2\eta-t_1,\dots,2\eta-t_m,-c)$.
Since generically $\psi_+$ and $\psi_-$ have different
multipliers, they are linearly independent over $K$: suppose
$\psi_+=C\psi_-$ for some $C\in K$. Then $C$ is $2\eta$-periodic
and obeys $C(\lambda+1)=e^{2c} C(\lambda)$. If $\eta$ is real and irrational
and $e^{2c}\neq 1$, we get a contradiction, so $\psi_+$ and $\psi_-$
are linearly independent. The case of generic $\eta$ is treated by
analytic continuation.

Next we show that every solution is a linear combination
of $\psi_+$ and $\psi_-$. The (difference) {\it Wronskian} $W(f,g)$ of two
functions $f$, $g$ is the function $f(\lambda+2\eta)g(\lambda)
-f(\lambda)g(\lambda+2\eta)$. Two meromorphic functions are
linearly dependent over $K$ if and only if their Wronskian 
vanishes. If $f,g$ are solutions of \Ref{lame} then their
Wronskian obeys the difference equation
$W(\lambda+2\eta)=u(\lambda)W(\lambda)$, where the function $u$ is
a combination of the coefficients of the \Ref{lame}.
Thus for any solution $f$,
the functions $A_\pm=\pm W(f,\psi_\mp)/W(\psi_+,\psi_-)$ are
$2\eta$-periodic.
On the other hand,   $A_\pm(\lambda)$ 
are (by Cramer's rule) the coefficients
in the expression  of $(f(\lambda),f(\lambda+2\eta))$ as 
 a linear combination of the linearly independent vectors
$(\psi_\pm(\lambda),\psi_\pm(\lambda+2\eta))$. In particular,
\be
f(\lambda)=A_+(\lambda)\psi_+(\lambda)+A_-(\lambda)\psi_-(\lambda).
\ee
Therefore $f$ is a linear combination of $\psi_\pm$ with 
$2\eta$-periodic coefficients.

It remains to prove that for each generic $\epsilon$ there are
not more than two solutions of the Bethe ansatz equation with given
$\epsilon$. Suppose that there were a third solution
$(t_1',\dots,t'_m,c')$ distinct from the two we have constructed.
In particular $(t'_1,\dots,t'_m)$ represents a point in the symmetric
power of the elliptic curve which is distinct from the two points
represented by 
$(t_1,\dots,t_m)$ and $(2\eta-t_1,\dots,2\eta-t_m)$. Let $\psi'$
be the corresponding eigenfunction. Thus $\psi'=a\psi_++b\psi_-$
with $2\eta$-periodic coefficients $a$, $b$. We have that
$a=W(\psi',\psi_-)/W(\psi_+,\psi_-)$. The three functions 
$\psi',\psi_+,\psi_-$ have all the form \Ref{psi3}. Therefore 
the meromorphic function $a$ obeys the equations 
$a(\lambda+2\eta)=a(\lambda)$,
$a(\lambda+1)=C_1 a(\lambda)$, $a(\lambda+\tau)=C_2 a(\lambda)$
for some constants $C_1,C_2$. If $\eta$ is generic, this implies that $a$ is of
the form $a_0e^{s\pi i/\eta}$ for some constant $a_0$ and
some integer $s$ ($a$ cannot have poles since it would have a dense set
of singularities). Similarly, $b$ has the same form. We may moreover assume
that $e^c$ is generic, so comparing multipliers we see that either
$a=0$ or $b=0$. But the zeros of $\psi'$ are not equal to the
zeros of $\psi_+$ or to the zeros of $\psi_-$,
a contradiction.
\end{prf}

We see from the proof of this result that the irreducible 
component(s) of $\bar X$ containing $P$ and
its reflected point $\sigma(P)$ form a {\em curve}
$Y\subset \bar X$. On this part of $\bar X$,
$\epsilon^2$ takes every generic value precisely twice. 

 Let us summarize our results.

\begin{thm} The closure $\bar X$ in $S^mE_\tau$ of the 
 Hermite--Bethe variety
contains an algebraic curve $Y$. It is 
a two-fold ramified covering
$\epsilon^2:Y\to {\Bbb P}^1$ of the Riemann sphere.
It has an involutive automorphism $\sigma:Y\to Y$
permuting the sheets and preserving $\epsilon^2$. For each
generic point  $t\in Y$, there are two solutions $(t,c)$,
$(t,c+i\pi/2\eta)$ of the Bethe ansatz equation in $Y\times
\C/\pi i\eta^{-1}\Z$. The corresponding eigenfunctions are related
by $\psi(\lambda)\to e^{\pi i\lambda/\eta}\psi(\lambda)$, and the
eigenvalues are the two square roots of $\epsilon^2(t)$.
\end{thm}

One way to formulate this result is that eigenfunctions
are parametrized by the {\em spectral curve}, the double covering
of $Y$ on which $\epsilon$ is single-valued.

\section{Proofs}\label{sproofs}

This section contains the proofs of Lemma \ref{sym} and
Theorem \ref{four}. These proofs are based on the following
technical result.

\begin{lemma}\label{tec} Let $V(z)$ with basis $e[{1}],e[{-1}]$
denote the fundamental representation
of $E_{\tau,\eta}(sl_2)$ with evaluation point $z$
(see Sect.\ \ref{firsts}). 
Let
$t_1,\dots,t_m$, $z_1,\dots,z_m$, $\tau,\eta$ be generic
complex numbers, such that $\operatorname{Im}(\tau)>0$, and let
$W=V(z_1)\otimes\cdots\otimes V(z_m)$.
Then the $2^m$ vectors
\be
\bigl(\prod_{j\in J}b(t_j)\bigr)(e[1]\otimes\cdots \otimes e[1])\in \Fun(W),
\ee
where $J$ runs over all subsets of $\{1,\dots, m\}$, are
linearly independent over the field of meromorphic functions
of $\lambda$.
\end{lemma}

\begin{prf}
Let us denote these vector by $w_J$. It is sufficient to show that
the values  $w_J(\lambda)$ are linearly independent for some value
of $\lambda$, and some value of the parameters. This will be shown
by considering the matrix relating $w_J(\lambda)$ to the basis
$e_J=e[{\sigma_1}]\otimes\cdots\otimes e[{\sigma_m}]$ with
$\sigma_j=1$ iff $j\in J$, and showing that, in some limit
of the parameters, this matrix is
upper triangular with respect to the lexicographical
ordering of binary numbers:
$J\geq K$ iff $\sum_{j\in J}2^j\geq\sum_{j\in K}2^j$.
The limit is obtained by first taking $\tau\to i\infty$ which
amounts to replacing $\theta(x)$ by $\sin(\pi x)$, then $\lambda\to
\infty$. In this limit, $b(t)$ acts on $V(z)$ as
$b(t)e[{-1}]=0$,
$b(t)e[1]=B(\eta)(1-\exp(-2\pi i(z-t+2\eta)))^{-1}e[{-1}]$ for some
constant $B(\eta)\neq 0$. 
Now let us choose $t_j=z_j+\epsilon$ for some fixed 
generic $\epsilon$, and set $z_j=\sqrt{-1}Nj$, $j=1,\dots,m$.
 Let $e[\sigma]_k$ be basis vectors of $V(z_k)$. Then,
as $N \to \infty$, the vector $b(t_j)e[1]_k$ tends to zero
if $j<k$ and tends to a nozero multiple of $e[-1]_k$ if $j=k$.
 On the other hand,
the matrix elements of $a$ and $d$  tend to
non-zero finite values in this limit. But the action
of $b(t)$ on a tensor product $u_1\otimes\cdots\otimes u_m$
consists of $m$ terms, the $j$'th one being given by the
action of $b$ on the $j$th factor and the action of $a$ or $d$
on the other $m-1$ factors (see \cite{FV}).
 
 It follows that $b(t_{j_1})\cdots b(t_{j_r})$ is a linear
combination of vectors $e_K$ where $K=\{k_1,\dots,k_r\}$
with $k_1\geq j_1,\dots,k_r\geq j_r$, and that the diagonal
matrix elements (such that $K=\{j_1,\dots,j_r\}$) are non-zero.
\end{prf}

\noindent{\it Proof of Lemma \ref{sym}}. The coefficients $A_j$, $D_j$ can
be in principle computed using the commutation relations
repeatedly, giving them as some universal polynomials in
the values of the functions $\alpha, \beta$, independent of
the highest weight modules we are considering. Since the $b$'s commute,
we obtain $m!$ ways of representing $a(w)b(t_1)\cdots b(t_m)$
 as a linear combination of 
$b(t_1)\cdots b(t_{j-1})b(w)b(t_{j+1})\cdots b(t_m)a(t_j)$
(and similarly for $d$).
Since by the previous Lemma there exists a vector in a module
so that these operators applied to this vector yield linearly
independent vectors, it follows that all these representations
must coincide. $\Box$

\vskip 30pt

\noindent{\it Proof of Theorem \ref{four}}. 
We have to write an explicit formula for
$b(t_1)\dots b(t_m)v_0$, where $v_0$ is a tensor product of highest
weight vectors of some highest weight modules. Let us first
consider the case where $v_0=u\otimes v$ is the tensor product
of two highest weight vectors. Then, using the rules for the
action on tensor products (see \cite{FV}) and the commutation
relations, we obtain a linear combination of terms of the
form
\be
\Gamma(-2\eta h^{(2)})[b(t_{i_1})\cdots b(t_{i_s})
a(t_{j_1})\cdots a(t_{j_{m-s}})]u \otimes 
 b(t_{k_1})...b(t_{k_{m-s}})
d(t_{l_1})...d(t_{l_{s}})v,
\ee
where $I=\{i_1 < ... < i_{s} \}$,
$J=\{ j_1 < ... < j_{m-s}\}$, 
$I \cup J = \{1,2,...,m\}$, 
$K=\{ k_1 < ... < k_{m-s}\}$, 
$L=\{ l_1 < ... < l_{s} \}$,
$K \cup L = \{1,2,...,m\}$. 
The notation $\Gamma(-2\eta h^{(2)})$ indicates that the
argument of the matrix valued function of $\lambda$ in the
square bracket must be shifted by $-2\eta$ times the weight
of the vector in the second factor.

The coefficients can be computed
using only the commutation relations and are therefore independent
of the choice of highest weight modules.
It is convenient to take $u$ and $v$ to be vectors of the form
of the previous Lemma, since it then follows from the linear
independence
that the coefficients are {\em uniquely} determined. 

Suppose that we compute the coefficients by first applying $b(t_m)$
to the tensor product, then $b(t_{m-1})$ and so on, and then use
the commutation relations to shift $d$'s and $a$'s to the right
of the $b$'s. Then it is clear that the coefficients of the
terms with $I\ni 1$ and $K\ni 1$ vanishes. Similarly, if 
we use the fact that the $b$'s commute to act with $b(t_j)$ at
the end, we see that the coefficients of the terms with
 $I\cap K\ni j$ must vanish. We conclude that the only
terms appearing with non-vanishing coefficients must have
$I\cap K=\emptyset$. Therefore
\be
b(t_1)\cdots b(t_m)(u\otimes v)
=\sum_I C_I\,\Gamma(-2\eta h^{(2)})\left[
\prod_{i\in I}b(t_i)\prod_{j\not\in I}a(t_j)\right]
u
\otimes \prod_{j\not\in I}b(t_j)\prod_{i\in I}d(t_i)v,
\ee
for some universal, uniquely defined,  coefficients $C_I(t_1,\dots,t_m)$,
$I\subset\{1,\dots,m\}$. By the commutativity of the $b$'s,
these coefficients can be computed in different ways with
the result that 
\be
C_{\sigma(I)}(t_1,\dots,t_m)=C_I(t_{\sigma(1)},\dots,t_{\sigma(m)}),
\ee
for all permutations $\sigma\in S_m$. Therefore, it is sufficient
to calculate $C_I$ for $I=\{1,\dots,s\}$, for all $s\in\{1,\dots, m\}$,
which can be done straightforwardly using the tensor product 
and commutation rules. Taking $u$ equal to the highest weight 
vector of an evaluation representation $V_\Lambda(z)$ and
$v$ a tensor product of such highest weight vectors gives a
recursive procedure to compute all coefficients given
in Theorem \ref{four}. $\Box$

\medskip
{\it Acknowledgement.} Part of this work was done at the
Erwin Schr\"odinger International Institute for Mathematical Physics.
We are grateful to the organizers of the program on ``Conformal and
integrable quantum field theory'', K. Gaw\c edzki and I. Todorov,
 for giving us this opportunity.

\end{document}